\documentclass{gtart}
\usepackage{amsmath,amssymb,epsf}
\input labelfig.tex
\input gtmonout
\volumenumber{1}
\volumeyear{1998}
\volumename{The Epstein birthday schrift}
\pagenumbers{167}{180}
\received{18 November 1997} 
\published{9 November 1998}
\papernumber{9}

\newtheorem{theorem}{Theorem}[section]
\newtheorem{lemma}[theorem]{Lemma}
\newtheorem{corollary}[theorem]{Corollary}

\newcommand{\la}{\alpha}
\newcommand{\lb}{\beta}
\newcommand{\ga}{\gamma}
\newcommand{\lG}{\Gamma}

\newcommand{\lv}{\varepsilon}

\newcommand{\lk}{\kappa}
\newcommand{\gl}{\lambda}

\newcommand{\lU}{\Upsilon}
\newcommand{\lO}{\Omega}

\newcommand{\h}{{\Bbb H}}

\newcommand{\R}{{\Bbb R}}
\newcommand{\Z}{{\Bbb Z}}

\newcommand{\cB}{{\cal B}}
\newcommand{\cC}{{\cal C}}
\newcommand{\cD}{{\cal D}}
\newcommand{\cE}{{\cal E}}
\newcommand{\cF}{{\cal F}}
\newcommand{\cH}{{\cal H}}
\newcommand{\cI}{{\cal I}}
\newcommand{\cJ}{{\cal J}}
\newcommand{\cK}{{\cal K}}
\newcommand{\cO}{{\cal O}}

\newcommand{\cN}{{\cal N}}

\newcommand{\cS}{{\cal S}}
\newcommand{\cT}{{\cal T}}
\newcommand{\sm}{\setminus}
\newcommand{\sub}{\subset}
\newcommand{\ra}{\rightarrow}
\newcommand{\Ra}{\Rightarrow}

\newcommand{\Lr}{\Leftrightarrow}

\newcommand{\ol}{\overline}

\newcounter{roman}

\newcounter{alpha}

\newcounter{arabic}
\newenvironment{arabiclist}{\begin{list}{{\rm(\arabic{arabic})}}{\usecounter{arabic}\setlength{\labelwidth}{\leftmargin}}}{\end{list}}

\begin{document}

\title[27 length inequalities define the Maskit domain in genus 2]{At most 27
length inequalities define Maskit's\\\vglue-2mm\\fundamental domain
for the modular group\\in genus 2} 
\covertitle{At most 27
length inequalities define Maskit's\\fundamental domain
for the modular group\\in genus 2} 
\asciititle{At most 27
length inequalities define Maskit's\\fundamental domain
for the modular group\\in genus 2} 

\author{David Griffiths} 

\address{Laboratoire de Mathematiques Pures de Bordeaux\\
Universite Bordeaux 1, 351 cours de la liberation\\
Talence 33405, Cedex, France}

\email{griffith@math.u-bordeaux.fr}


\begin{abstract}

In recently published work Maskit constructs a fundamental domain
$\cD_g$ for the Teichm\"{u}ller modular group of a closed surface
$\cS$ of genus $g\geq 2$.  Maskit's technique is to demand that a
certain set of $2g$ non-dividing geodesics $\cC_{2g}$ on $\cS$
satisfies certain shortness criteria. This gives an a priori infinite
set of length inequalities that the geodesics in $\cC_{2g}$ must
satisfy. Maskit shows that this set of inequalities is finite and that
for genus $g=2$ there are at most 45. In this paper we improve this
number to 27. Each of these inequalities: compares distances between
Weierstrass points in the fundamental domain $\cS\sm\cC_4$ for
$\cS$; and is realised (as an equality) on one or other of two special
surfaces.

\end{abstract} 
\asciiabstract{%
In recently published work Maskit constructs a fundamental domain
D_g for the Teichmueller modular group of a closed surface
S of genus g>1.  Maskit's technique is to demand that a
certain set of 2g non-dividing geodesics C_{2g} on S
satisfies certain shortness criteria. This gives an a priori infinite
set of length inequalities that the geodesics in C_{2g} must
satisfy. Maskit shows that this set of inequalities is finite and that
for genus g=2 there are at most 45. In this paper we improve this
number to 27. Each of these inequalities: compares distances between
Weierstrass points in the fundamental domain S-C_4 for
S; and is realised (as an equality) on one or other of two special
surfaces.}

\maketitle
\setcounter{section}{-1} \section{Introduction and preliminaries}

In this paper we consider a fundamental domain defined by Maskit in
\cite{Mas3} for the action of the Teichm\"{u}ller modular group on the
Teichm\"{u}ller space of a closed surface of genus $g\geq 2$ in the
special case of genus $g=2$. McCarthy and Papadopoulos~\cite{M-P} have
also defined such a fundamental domain, modelled on a Dirichlet region; for
punctured surfaces there is the celebrated cell decomposition and
associated fundamental domain due to Penner~\cite{Pen}. For genus
$g=2$ Semmler~\cite{Sem} has defined a fundamental domain based on
locating the shortest dividing geodesic. Also for low signature
surfaces the reader is referred to the papers of Keen~\cite{Keen} and
of Maskit \cite{Mas2},~\cite{Mas3}.

Throughout $\cS$ will denote a closed orientable surface of genus
$g=2$, with some fixed hyperbolic metric. We say that a simple closed
geodesic $\ga$ on $\cS$ is: {\sl dividing\/} if $\cS\setminus\ga$ has two
components; or {\sl non-dividing\/} if $\cS\setminus\ga$ has one component. By
{\em non-dividing geodesic} we shall always mean simple closed
non-dividing geodesic. We denote the length of $\ga$ with respect to
the hyperbolic metric on $\cS$ by $l(\ga)$. Let $|\la\cap\lb|$ denote
the number of intersection points of two distinct geodesics $\la,\lb$.

We define a {\em chain} $\cC_n=\ga_1,\ldots,\ga_n$ to be an ordered
set of non-dividing geodesics such that: $|\ga_i\cap\ga_{i+1}|=1$ for
$1\leq i\leq n-1$ and $\ga_i\cap\ga_{i'}=\emptyset$ otherwise. We say
that a chain $\cC_n$ has {\em length} $n$, where $1\leq n\leq
5$. Likewise we define a {\em bracelet} $\cB_n=\ga_1,\ldots,\ga_n$ to
be an ordered set of non-dividing geodesics such that:
$|\ga_i\cap\ga_{i+1}|=1$ for $1\leq i\leq n-1,~|\ga_n\cap\ga_1|=1$ and
$\ga_i\cap\ga_{i'}=\emptyset$ otherwise. Again we say that $\cB_n$ has
{\em length} $n$, where $3\leq n\leq 6$. Following Maskit, we call a
bracelet of length 6 a {\em necklace}.

For $n\leq 4$ a chain of length $n$ can be always be extended to a
chain of length $n+1$. For $n=4$ this extension is unique. Likewise a
chain of length 5 extends uniquely to a necklace. So chains of length
4 or 5 and necklaces can be considered equivalent. We shall usually
work with length 4 chains, which we call {\em standard}. (Maskit, for
genus $g$, usually works with chains of length 2g+1, which he calls
standard.)

As Maskit shows in \cite{Mas3} each surface, standard chain pair
$\cS,\cC_4$ gives a canonical choice of generators for the Fuchsian
group $F$ such that $\h^2/F=\cS$ and hence a point in
$\cD\cF(\pi_1(\cS),PSL(2,\R))$, the set of discrete faithful
representations of $\pi_1(\cS)$ into $PSL(2,\R)$. Essentially this
representation corresponds to the fundamental domain $\cS\sm\cC_4$
together with orientations for its side pairing elements. As Maskit
observes, it is well known that $\cD\cF(\pi_1(\cS),PSL(2,\R))$ is real
analytically equivalent to Teichm\"{u}ller space. So, we define the
{\em Teichm\"{u}ller space} of closed orientable genus $g=2$ surfaces
$\cT_2$ to be the set of pairs $\cS,\cC_4$.

We say that a standard chain $\cC_4=\ga_1,\ldots,\ga_4$ is {\em
minimal} if for any chain $\cC'_m=\ga_1,\ldots,\ga_{m-1},\la_m$ we
have $l(\ga_m)\leq l(\la_m)$ for $1\leq m\leq 4$. We then define the
{\em Maskit domain} $\cD_2\sub\cT_2$ to be the set of surface,
standard chain pairs $\cS,\cC_4$ with $\cC_4$ minimal.

For $\cC_4$ to be minimal the geodesics $\ga_1,\ldots,\ga_4$ must
satisfy an a priori infinite set of length inequalities. For genus
$g$, Maskit gives an algorithm using cut-and-paste to show that only a
finite number $N_g$ of length inequalities need to be
satisfied. Applying his algorithm to genus $g=2$, Maskit showed that
$N_2\leq 45$. We establish an independent proof that $N_2\leq 27$. We
could have shown that 18 of Maskit's 45 inequalities follow from the
other 27. However, by tayloring all our techniques to the special case of
genus 2, we are able to produce a much shorter proof.

The fact that 18 of Maskit's 45 inequalities follow from the other 27
follows from applications of Theorem~\ref{tri} (which appeared as
Theorem~1.1 in~\cite{Gri1}) and of Corollary~\ref{angle
corollary}. The latter follows immediately from Theorem~\ref{angle
characterisation}, for which we give a proof in this paper. This is a
characterisation of the octahedral surface $\cO ct$ (the well known
genus two surface of maximal symmetry group) in terms of a finite set
of length inequalities.

The 27 length inequalities have the properties that: each is realised
on one or other of two special surfaces (for all but $2$ this special
surface is $\cO ct$); and each compares distances between
Weierstrass points in the fundamental domain $\cS\sm\cC_4$ for
$\cS$.

The author would like to thank Bill Harvey, Bernie Maskit, Peter
Buser, Klaus-Dieter Semmler and Christophe Bavard for hospitality and
helpful discussions.  The author was supported for this work by the
Swiss National Science Foundation on a Royal Society Exchange
Fellowship at EPFL, Lausanne, Switzerland and is currently supported
by the French Government as a boursier on a Sejour Scientifique.

\section{The hyperelliptic involution and the main result}

It is well known that every closed genus two surface without boundary
$\cS$ admits a uniquely determined {\em hyperelliptic involution}, an
isometry of order two with six fixed points, which we denote by
$\cJ$. The fixed points of $\cJ$ are known as {\em Weierstrass
points}. Every simple closed geodesic $\ga\sub\cS$ is setwise fixed by
$\cJ$, and the restriction of $\cJ$ to $\ga$ has no fixed points if
$\ga$ is dividing and two fixed points if $\ga$ is non-dividing (see
Haas--Susskind~\cite{Haas-Sus}). So every non-dividing geodesic on
$\cS$ passes through two Weierstrass points. It is a simple
consequence that sequential geodesics in a chain intersect at
Weierstrass points. We say that two non-dividing geodesics
$\la,\lb$ {\em cross} if $\la\neq\lb$ and $\la\cap\lb$ contains a
point that is not a Weierstrass point.

The quotient {\em orbifold} $\cO\cong\cS/\cJ$ is a sphere with six
order two cone points, endowed with a fixed hyperbolic metric. Each
cone point on $\cO$ is the image of a Weierstrass point under the
projection $\cJ\co\cS\ra\cO$ and each non-dividing geodesic on $\cS$
projects to a simple geodesic between distinct cone points on $\cO$ --
what we shall call an {\em arc}. Definitions of chains, bracelets and
crossing all pass naturally to the quotient.

Let $\cC_4$ be a standard chain on $\cS$, which extends to a necklace
$\cN$. We number Weierstrass points on $\cN$ so that
$\omega_i=\ga_{i-1}\cap\ga_i$ for $2\leq i\leq 6$ and
$\omega_1=\ga_6\cap\ga_1$. Choose an orientation upon $\cS$ and
project to the quotient orbifold $\cO\cong\cS/\cJ$ -- for the rest of
the paper we shall work on the quotient orbifold $\cO$. We label the
components of $\cO\sm\cN$ by $\cH,\ol{\cH}$ so that
$\ga_1,\ldots,\ga_6$ lie anticlockwise around $\cH$. Label by
$\lb_{j,k}^{i_1,i_2,\ldots,i_n}$
(respectively $\ol{\lb_{j,k}^{i_1,i_2,\ldots,i_n}}$) the arc between the cone
points $\omega_j,\omega_k$ ($j<k$) crossing the sequence of arcs
$\ga_{i_1},\ga_{i_2},\ldots,\ga_{i_n}$ and having the subarc between
$\omega_j,\ga_{i_1}$ lying in $\cH$ (respectively $\ol{\cH}$).

Our main result is then the following. (We abuse notation so that
$\lb_{1,6}=\ol{\lb_{1,6}}=\ga_6$ and
$\lb_{2,3}=\ol{\lb_{2,3}}=\ga_2$. We then have repetitions,
$l(\ga_2)\leq l(\ga_6)$ twice, and redundancies, $l(\ga_2)\leq
l(\ga_2)$ also twice.)

\begin{theorem}\label{probably}

The standard chain $\cC_4$ is minimal if the following are satisfied:

\begin{arabiclist}

\item $l(\ga_1)\leq l(\ga_i),i\in\{2,3,4,5\}$

\item $l(\ga_2)\leq l(\lb_{i,j}),l(\ol{\lb_{i,j}}),
l(\lb_{2,5}^6),l(\ol{\lb_{2,5}^6})$, $i\in\{1,2\}$, $j\in\{3,4,5,6\}$

\item $l(\ga_3)\leq l(\lb_{3,j}),l(\ol{\lb_{3,j}}),l(\lb_{3,4}^6),
l(\ol{\lb_{3,4}^6})$, $j\in\{5,6\}$

\item $l(\ga_4)\leq l(\lb_{4,6}),l(\ol{\lb_{4,6}})$.

\end{arabiclist}

\end{theorem}

Each length $l(\ga_i)$ or $l(\lb_{j,k})$ (respectively $l(\ol{\lb_{j,k}})$)
is a distance between cone points in $\cH$ (respectively $\ol{\cH}$). Likewise
each length $l(\lb^6_{j,k})$, $l(\ol{\lb^6_{j,k}})$ is a distance between
cone points in $\cO\sm\cC_5$. So each length inequality in
Theorem~\ref{probably} compares distances between cone points in
$\cO\sm\cC_5$ (and hence distances between Weierstrass points in 
$\cS\sm\cC_4$).
 
\begin{figure}[htbp] 
\cl{%
\SetLabels 
(-0.005*0.48)\mbox{\scriptsize $\omega_5$}\\
(0.114*0.44)\mbox{\scriptsize $\omega_6$}\\
(0.235*0.38)\mbox{\scriptsize $\omega_3$}\\
(0.13*1.02)\mbox{\scriptsize $\omega_4$}\\
(0.13*-0.07)\mbox{\scriptsize $\omega_2$}\\
(0.1375*0.59)\mbox{\scriptsize $\omega_1$}\\
(0.37*0.44)\mbox{\scriptsize $\omega_6$}\\
(0.495*0.42)\mbox{\scriptsize $\omega_3$}\\
(0.3775*1.02)\mbox{\scriptsize $\omega_4$}\\
(0.3775*-0.07)\mbox{\scriptsize $\omega_1$}\\
(0.39*0.59)\mbox{\scriptsize $\omega_2$}\\
(0.265*0.56)\mbox{\scriptsize $\omega_5$}\\
(0.6225*0.44)\mbox{\scriptsize $\omega_5$}\\
(0.7475*0.41)\mbox{\scriptsize $\omega_2$}\\
(0.635*1.02)\mbox{\scriptsize $\omega_4$}\\
(0.635*-0.07)\mbox{\scriptsize $\omega_1$}\\
(0.52*0.56)\mbox{\scriptsize $\omega_6$}\\
(0.6425*0.59)\mbox{\scriptsize $\omega_3$}\\
(0.94*0.47)\mbox{\scriptsize $\omega_6$}\\
(1.015*0.46)\mbox{\scriptsize $\omega_5$}\\
(0.9*0.02)\mbox{\scriptsize $\omega_3$}\\
(0.9*0.93)\mbox{\scriptsize $\omega_4$}\\
(0.7775*0.56)\mbox{\scriptsize $\omega_1$}\\
(0.8225*0.48)\mbox{\scriptsize $\omega_2$}\\
\endSetLabels 
\AffixLabels{\hbox to 4.74truein{\hfil\epsfxsize 4.7truein
\epsffile{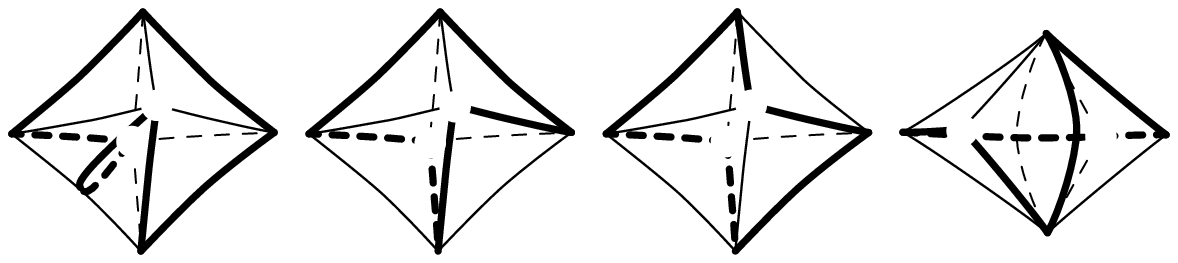}}}}
\caption{How the length inequalities in Theorem~\ref{probably} are
realized on $\cO ct$ and $\cE$}
\label{Fig:realise} 
\end{figure}

Theorem~\ref{probably} gives a sufficient list of inequalities. As to
the necessity each inequality, we make the following observation. Each
inequality is realised (as an equality) on either $\cO ct$ or $\cE$ -- %
cf Theorem~1.1 in~\cite{Gri2}. The octahedral orbifold $\cO ct$ is
the well known orbifold of maximal conformal symmetry group. Any
minimal standard chain on $\cO ct$ lies in its set of shortest
arcs. This arc set has the combinatorial edge pattern of the Platonic
solid. The exceptional orbifold $\cE$, which was constructed in
\cite{Gri2}, has conformal symmetry group $\Z_2\times\Z_2$. However it
is not defined by the action of its symmetry group alone, it also
requires a certain length inequality to be satisfied. Any minimal
standard chain on $\cE$ lies in its set of shortest and second
shortest arcs.

In Figure~\ref{Fig:realise} we have illustrated necklaces on $\cO ct$
and $\cE$ that are the extentions of minimal standard chains. As with
other figures in this paper, we use wire frame diagrams to illustrate
the orbifolds. Solid (respectively dashed) lines represent arcs in front
(respectively behind) the figure. Thick lines represent arcs in the necklace
$\cN$. The minimal standard chain on $\cE$ in Figure~\ref{Fig:realise}
has: $l(\ga_1)=l(\ga_5)$;
$l(\ga_2)=l(\ol{\lb_{1,3}})=l(\lb_{1,4})=l(\lb_{2,4})$;
$l(\ga_3)=l(\ol{\lb_{3,4}^6})=l(\ol{\lb_{3,5}})=l(\ol{\lb_{3,6}})$;
$l(\ga_4)=l(\lb_{4,6})$. Making such a list for all the orbifolds in
Figure~\ref{Fig:realise}, together with their mirror images, we see that
all the inequalities in Therem~\ref{probably} are realised as
equalities on either $\cO ct$ or $\cE$.

\section{Length inequalities for systems of arcs}\label{systems}

In order to prove Theorem~\ref{probably} we need a number of length
inequality results for systems of arcs. Let
$\cK_4=\lk_{0,1},\ldots,\lk_{3,0}$ denote a length 4 bracelet such
that each component of $\cO\sm\cK_4$ contains an interior cone
point. Using mod 4 addition throughout, label cone points: on $\cK_4$
by $c_k=\lk_{k-1,k}\cap\lk_{k,k+1}$ for $k\in\{0,\ldots,3\}$; and off
$\cK_4$ by $c_l$ for $l\in\{4,5\}$. Label by $\cO_l$ the component of
$\cO\sm\lU$ containing $c_l$ and label arcs in $\cO_l$ so that
$\lk_{k,l}$ is between $c_k,c_l$. Let $\gl_k$ denote the arc between
$c_4,c_5$ crossing only $\lk_{k,k+1}\sub\cK_4$.

The following two results appeared as Lemma~2.3 in~\cite{Gri2} (in
Maskit's terminology this is a cut-and-paste) and as Theorem~1.1
in~\cite{Gri1} respectively.

\begin{lemma}\label{c&p-2}  

{\rm (i)}\qua $2l(\lk_{0,4})<l(\gl_0)+l(\gl_3)$\quad {\rm (ii)}
\qua $2l(\lk_{3,0})<l(\gl_0)+l(\gl_2)$.

\end{lemma}

\begin{theorem}\label{tri} 

If $l(\lk_{3,4})\leq l(\lk_{0,4})$, $l(\lk_{3,5})\leq l(\lk_{0,5})$,
$l(\gl_0)\leq l(\gl_2)$ then\break $l(\lk_{3,4})=l(\lk_{0,4})$,
$l(\lk_{3,5})=l(\lk_{0,5})$, $l(\gl_0)=l(\gl_2)$.

\end{theorem}    

\begin{corollary}\label{obvious} 

If $l(\lk_{3,4})\leq l(\lk_{0,4})$, $l(\lk_{3,5})\leq l(\lk_{0,5})$,
$l(\lk_{1,4})\leq l(\lk_{2,4})$ then $l(\lk_{1,5})\geq l(\lk_{2,5}).$

\end{corollary}

\begin{proof}[Proof of Corollary~\ref{obvious}] Since 
$l(\lk_{3,4})\leq l(\lk_{0,4})$, $l(\lk_{3,5})\leq l(\lk_{0,5})$
Theorem~\ref{tri} implies that $l(\gl_0)\geq l(\gl_2)$. Moreover
$l(\lk_{1,4})\leq l(\lk_{2,4})$ and so again, by Theorem~\ref{tri},
$l(\lk_{1,5})\geq l(\lk_{2,5})$.\end{proof}

\begin{figure}[htbp] 
\begin{center} 
\leavevmode 
\SetLabels 
(0.1325*0.6)\mbox{\scriptsize $c_0$}\\
(0.235*0.43)\mbox{\scriptsize $c_1$}\\
(0.11375*0.44)\mbox{\scriptsize $c_2$}\\
(0.0*0.48)\mbox{\scriptsize $c_3$}\\
(0.125*-0.07)\mbox{\scriptsize $c_4$}\\
(0.125*1.02)\mbox{\scriptsize $c_5$}\\
(0.3875*0.59)\mbox{\scriptsize $c_0$}\\
(0.495*0.42)\mbox{\scriptsize $c_1$}\\
(0.3675*0.44)\mbox{\scriptsize $c_2$}\\
(0.265*0.53)\mbox{\scriptsize $c_3$}\\
(0.3775*-0.07)\mbox{\scriptsize $c_4$}\\
(0.3775*1.02)\mbox{\scriptsize $c_5$}\\
(0.6425*0.59)\mbox{\scriptsize $c_0$}\\
(0.7475*0.42)\mbox{\scriptsize $c_1$}\\
(0.625*0.44)\mbox{\scriptsize $c_2$}\\
(0.52*0.53)\mbox{\scriptsize $c_3$}\\
(0.635*-0.07)\mbox{\scriptsize $c_4$}\\
(0.635*1.02)\mbox{\scriptsize $c_5$}\\
(0.9*0.59)\mbox{\scriptsize $c_0$}\\
(1.015*0.46)\mbox{\scriptsize $c_1$}\\
(0.88*0.44)\mbox{\scriptsize $c_2$}\\
(0.7775*0.55)\mbox{\scriptsize $c_3$}\\
(0.885*-0.06)\mbox{\scriptsize $c_4$}\\
(0.885*1.02)\mbox{\scriptsize $c_5$}\\
\endSetLabels 
\AffixLabels{ \epsffile{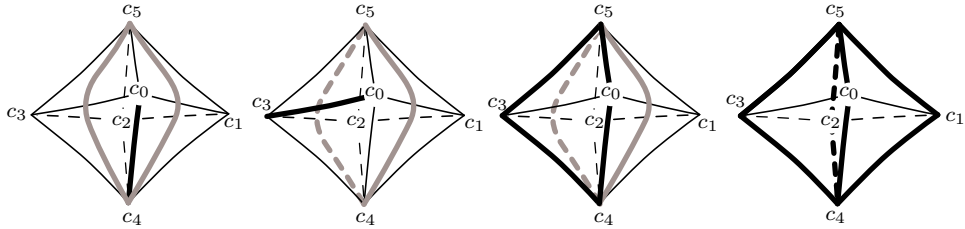}} 
\end{center}
\caption{Arc sets for Lemma~\ref{c&p-2}, for Theorem~\ref{tri} and for
Corollary~\ref{obvious}}
\label{Fig:length_2_3} 
\end{figure}

\begin{theorem}\label{angle characterisation} 

Suppose $l(\lk_{2,3})\leq l(\lk_{2,l})$, $l(\lk_{3,0})\leq  
l(\lk_{1,2})\leq
\{l(\lk_{0,l}),l(\lk_{1,l})\}$ and $l(\lk_{0,1})\leq
\{l(\lk_{0,l}),l(\lk_{3,l})\}$ then $l(\lk_{k,l})=l(\lk_{k,k+1})$ for each
$k,l$ and $\cO$ is the octahedral orbifold.

\end{theorem} 

\begin{proof}[Proof of Theorem~\ref{angle characterisation}]
 We postpone this until
Section~\ref{proofs}.\end{proof}

\begin{corollary}\label{angle corollary} 

Suppose $l(\lk_{2,3})\leq l(\lk_{2,l})$, $l(\lk_{1,2})\leq
\{l(\lk_{0,l}),l(\lk_{1,l})\}$ and\break $l(\lk_{0,1})\leq
\{l(\lk_{0,l}),l(\lk_{3,l})\}$ then $l(\lk_{3,0})\geq l(\lk_{1,2})$.

\end{corollary}

\begin{proof}[Proof of Corollary~\ref{angle corollary}] If $l(\lk_{3,0})\leq
l(\lk_{1,2})$ then by Theorem~\ref{angle characterisation}
$l(\lk_{k,l})=l(\lk_{k,k+1})$ for each $k,l$. In particular
$l(\lk_{3,0})=l(\lk_{1,2})$. So $l(\lk_{3,0})\geq l(\lk_{1,2})$.\end{proof}

\section{The proofs}\label{proofs}

\begin{proof}[Proof of Theorem~\ref{probably}] Let $\la_m$ denote an arc such
that $\cC'_m=\ga_1,\ldots,\ga_{m-1},$ $\la_m$ is a chain, for
$1\leq m\leq 4,~\la_m\neq\ga_m$. We will show that $l(\ga_m)\leq
l(\la_m)$ for arcs of the form $\lb_{j,k}^{i_1,i_2,\ldots,i_n}$. The
same arguments work for arcs of the form
$\ol{\lb_{j,k}^{i_1,i_2,\ldots,i_n}}$. Let $X(\la,\lb)$ denote the
number of crossing points of a distinct pair of arcs $\la,\lb$ -- %
ie the number of intersection points of $\la,\lb$ that are not cone
points. Let $n=\infty$, if $X(\ga_m,\la_i)=0$ for
$i\in\{1,\ldots,6\}$; otherwise, let $n=~{\rm min}~i\in\{1,\ldots,6\}$ such
that $X(\ga_n,\la_i)>0$. We note that $n\geq m$.

Let $P_{m,n,p}$ be the proposition that $l(\ga_m)\leq
l(\la_m)$ for $X(\la_m,\ga_n)=p$. Clearly, if $n=\infty$
then $p=0$. For $n\in\{5,6\}$ it is not hard to show that $p=1$. For
$n\in\{1,\ldots,4\}$ we consider $p=1$ and $p>1$. We order the
propositions as follows: $P_{4,\infty,0},\ldots,P_{1,\infty,0}$ which
is followed by $P_{4,6,1},P_{4,5,1},\ldots,P_{1,6,1},P_{1,5,1}$
followed by $P_{4,4,1},P_{4,4,p>1}$ which is followed by
$P_{3,4,1},P_{3,4,p>1},P_{3,3,1},P_{3,3,p>1}$ followed by $
P_{2,4,1},P_{2,4,p>1},\ldots,P_{2,2,1},P_{2,2,p>1}$ followed by
$P_{1,4,1},P_{1,4,p>1},\ldots,P_{1,1,1},P_{1,1,p>1}$.

Suppose $n=\infty,~\la_m$ does not cross $\cN$. If $m>1$ then
$P_{m,\infty,0}$ is a hypothesis. If $m=1$ then either
$P_{1,\infty,0}$ is a hypothesis, $\la_1=\ga_i$ for some
$i\in\{2,\ldots,5\}$, or $P_{1,\infty,0}$ follows from the hypotheses,
$l(\ga_1)\leq l(\ga_i),l(\ga_i)\leq l(\la_1)$ for some
$i\in\{2,3,4\}$.

Suppose $n\in\{5,6\},~\la_m$ crosses $\cN$ but does not cross
$\cC_4$. 

For $m=4$, by inspection, $\la_4=\lb_{4,5}^6$. So $\la_m,\ga_m$ share
endpoints, $n>m+1$ and we can apply the argument (i) below. So we have
$P_{4,n,1}$ for $n\in\{5,6\}$.

In Figures 3,4,5 we illustrate applications of
length inequalities results to the proof. As above we use wire frame
figures of the octahedral orbifold, with the necklace $\cN$ in thick
black. Other arcs are in thick grey. Figures have been drawn so arcs
in the application correspond to arcs in the length inequality result.

\begin{figure}[htbp] 
\begin{center} 
\leavevmode 
\SetLabels 
(0.0*0.48)\mbox{\scriptsize $\omega_1$}\\
(0.11375*0.44)\mbox{\scriptsize $\omega_2$}\\
(0.235*0.43)\mbox{\scriptsize $\omega_3$}\\
(0.125*1.02)\mbox{\scriptsize $\omega_4$}\\
(0.125*-0.07)\mbox{\scriptsize $\omega_5$}\\
(0.1325*0.6)\mbox{\scriptsize $\omega_6$}\\
(0.3675*0.44)\mbox{\scriptsize $\omega_1$}\\
(0.495*0.42)\mbox{\scriptsize $\omega_2$}\\
(0.3775*1.02)\mbox{\scriptsize $\omega_3$}\\
(0.3775*-0.07)\mbox{\scriptsize $\omega_4$}\\
(0.3875*0.59)\mbox{\scriptsize $\omega_5$}\\
(0.265*0.53)\mbox{\scriptsize $\omega_6$}\\
(0.625*0.44)\mbox{\scriptsize $\omega_1$}\\
(0.7475*0.42)\mbox{\scriptsize $\omega_2$}\\
(0.635*1.02)\mbox{\scriptsize $\omega_3$}\\
(0.635*-0.07)\mbox{\scriptsize $\omega_4$}\\
(0.52*0.53)\mbox{\scriptsize $\omega_5$}\\
(0.6425*0.59)\mbox{\scriptsize $\omega_6$}\\
(0.88*0.44)\mbox{\scriptsize $\omega_1$}\\
(1.015*0.46)\mbox{\scriptsize $\omega_2$}\\
(0.885*-0.06)\mbox{\scriptsize $\omega_3$}\\
(0.885*1.02)\mbox{\scriptsize $\omega_4$}\\
(0.7775*0.55)\mbox{\scriptsize $\omega_5$}\\
(0.9*0.59)\mbox{\scriptsize $\omega_6$}\\
\endSetLabels 
\AffixLabels{ \epsffile{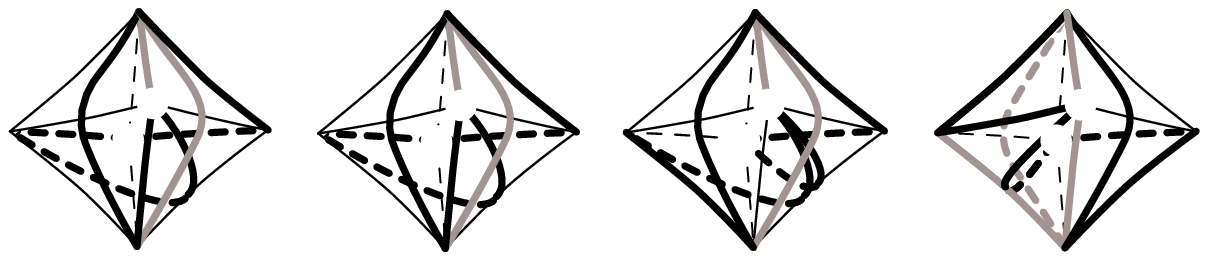}} 
\end{center}
\caption{Application (i) for
$\la_4=\lb_{4,5}^6,\la_3=\lb_{3,4}^5$ and $\lb_{3,4}^{6,5}$ and of
Theorem~\ref{tri}, (ii) for $\la_3=\lb_{3,5}^6$}
\label{Fig:arcs_4,3} 
\end{figure}

For $m=3$. By inspection, $\la_3$ is one of
$\lb_{3,4}^5,\lb_{3,4}^6,\lb_{3,4}^{6,5},\lb_{3,5}^6$. For
$\lb_{3,4}^5,\lb_{3,4}^{6,5},\lb_{3,4}^6$: $\ga_m,\la_m$ share
endpoints, $n>m+1$ and so we can apply either argument (i) or (ii)
below. For $\lb_{3,5}^6$ we can apply Theorem~\ref{tri} in conjunction
with argument (ii): by hypothesis $l(\ga_4)\leq l(\lb_{4,6})$ and by
argument (ii) $l(\ga_3)\leq l(\lb_{3,4}^6)$ and so $l(\lb_{3,5}^6)\geq
l(\lb_{3,6})$. Again by hypothesis $l(\ga_3)\leq l(\lb_{3,6})$ and so
$l(\ga_3)\leq l(\lb_{3,6})\leq l(\lb_{3,5}^6)$. This gives $P_{3,n,1}$
for $n\in\{5,6\}$. 

For $m=2,~\la_2$ is one of $\lb_{2,3}^5,\lb_{2,3}^6,\lb_{2,3}^{6,5}$
or one of
$\lb_{2,4}^5,\lb_{2,4}^6,\lb_{2,4}^{6,5},\lb_{2,5}^6,\lb_{1,3}^5,
\lb_{1,4}^5$. By hypothesis $l(\ga_2)\leq l(\lb_{2,5}^6)$. For
$\lb_{2,3}^6,\lb_{2,3}^{6,5},\lb_{2,3}^5$ we can again apply either
argument (i) or (ii). For
$\lb_{2,4}^5,\lb_{2,4}^6,\lb_{2,4}^{6,5},\lb_{1,3}^5$ we apply
Theorem~\ref{tri} in conjunction with argument (ii). We give the
argument for $\lb_{2,4}^5$. By argument (ii), we have
$l(\ga_2)<l(\lb_{2,3}^5)$. Also, by hypothesis, $l(\ga_3)\leq
l(\lb_{3,5})$ and so by Theorem~\ref{tri}
$l(\lb_{2,5})<l(\lb_{2,4}^5)$.  Again, by hypothesis, $l(\ga_2)\leq
l(\lb_{2,5})$ and so $l(\ga_2)\leq l(\lb_{2,5})<l(\lb_{2,4}^5)$.

For $\la_2=\lb_{1,4}^5$ we argue as follows. By hypothesis we have
$l(\ga_3)\leq l(\lb_{3,5}),l(\ol{\lb_{3,6}})$ and $l(\ga_2)\leq
l(\lb_{1,5}),l(\ga_6),l(\lb_{2,5}),l(\ol{\lb_{2,6}})$ and
$l(\ga_1)\leq l(\lb_{1,5}),l(\ga_6),l(\ga_4),l(\ol{\lb_{4,6}})$. By
Corollary~\ref{angle corollary}: $l(\lb_{1,4}^5)\geq l(\ga_2)$. Hence
$P_{2,n,1}$ for $n\in\{5,6\}$.

\begin{figure}[htbp] 
\begin{center} 
\leavevmode 
\SetLabels 
(0.23*0.735)\mbox{\scriptsize $\omega_1$}\\
(0.12*0.535)\mbox{\scriptsize $\omega_2$}\\
(0.12*1.01)\mbox{\scriptsize $\omega_3$}\\
(-0.01*0.775)\mbox{\scriptsize $\omega_4$}\\
(0.1275*0.82)\mbox{\scriptsize $\omega_5$}\\
(0.1075*0.7575)\mbox{\scriptsize $\omega_6$}\\
(0.385*0.82)\mbox{\scriptsize $\omega_1$}\\
(0.3775*1.01)\mbox{\scriptsize $\omega_2$}\\
(0.3775*0.535)\mbox{\scriptsize $\omega_3$}\\
(0.4875*0.735)\mbox{\scriptsize $\omega_4$}\\
(0.365*0.7575)\mbox{\scriptsize $\omega_5$}\\
(0.26*0.8)\mbox{\scriptsize $\omega_6$}\\
(0.745*0.74)\mbox{\scriptsize $\omega_1$}\\
(0.635*1.01)\mbox{\scriptsize $\omega_2$}\\
(0.635*0.535)\mbox{\scriptsize $\omega_3$}\\
(0.52*0.795)\mbox{\scriptsize $\omega_4$}\\
(0.6275*0.7575)\mbox{\scriptsize $\omega_5$}\\
(0.64875*0.82)\mbox{\scriptsize $\omega_6$}\\
(1.0125*0.775)\mbox{\scriptsize $\omega_1$}\\
(0.8875*0.535)\mbox{\scriptsize $\omega_2$}\\
(0.8875*1.01)\mbox{\scriptsize $\omega_3$}\\
(0.7725*0.81)\mbox{\scriptsize $\omega_4$}\\
(0.8975*0.82)\mbox{\scriptsize $\omega_5$}\\
(0.8775*0.7575)\mbox{\scriptsize $\omega_6$}\\
(0.1275*0.255)\mbox{\scriptsize $\omega_1$}\\
(0.12*-0.03)\mbox{\scriptsize $\omega_2$}\\
(0.12*0.445)\mbox{\scriptsize $\omega_3$}\\
(-0.01*0.20)\mbox{\scriptsize $\omega_4$}\\
(0.11*0.195)\mbox{\scriptsize $\omega_5$}\\
(0.23*0.17)\mbox{\scriptsize $\omega_6$}\\
(0.4875*0.17)\mbox{\scriptsize $\omega_1$}\\
(0.3775*-0.03)\mbox{\scriptsize $\omega_2$}\\
(0.3775*0.445)\mbox{\scriptsize $\omega_3$}\\
(0.26*0.235)\mbox{\scriptsize $\omega_4$}\\
(0.3675*0.195)\mbox{\scriptsize $\omega_5$}\\
(0.385*0.255)\mbox{\scriptsize $\omega_6$}\\
(0.635*0.445)\mbox{\scriptsize $\omega_1$}\\
(0.635*-0.03)\mbox{\scriptsize $\omega_2$}\\
(0.52*0.24)\mbox{\scriptsize $\omega_3$}\\
(0.625*0.195)\mbox{\scriptsize $\omega_4$}\\
(0.745*0.175)\mbox{\scriptsize $\omega_5$}\\
(0.645*0.255)\mbox{\scriptsize $\omega_6$}\\
(0.8875*-0.03)\mbox{\scriptsize $\omega_1$}\\
(1.0125*0.2)\mbox{\scriptsize $\omega_2$}\\
(0.8875*0.445)\mbox{\scriptsize $\omega_3$}\\
(0.7725*0.24)\mbox{\scriptsize $\omega_4$}\\
(0.9*0.255)\mbox{\scriptsize $\omega_5$}\\
(0.8775*0.195)\mbox{\scriptsize $\omega_6$}\\
\endSetLabels 
\AffixLabels{\epsffile{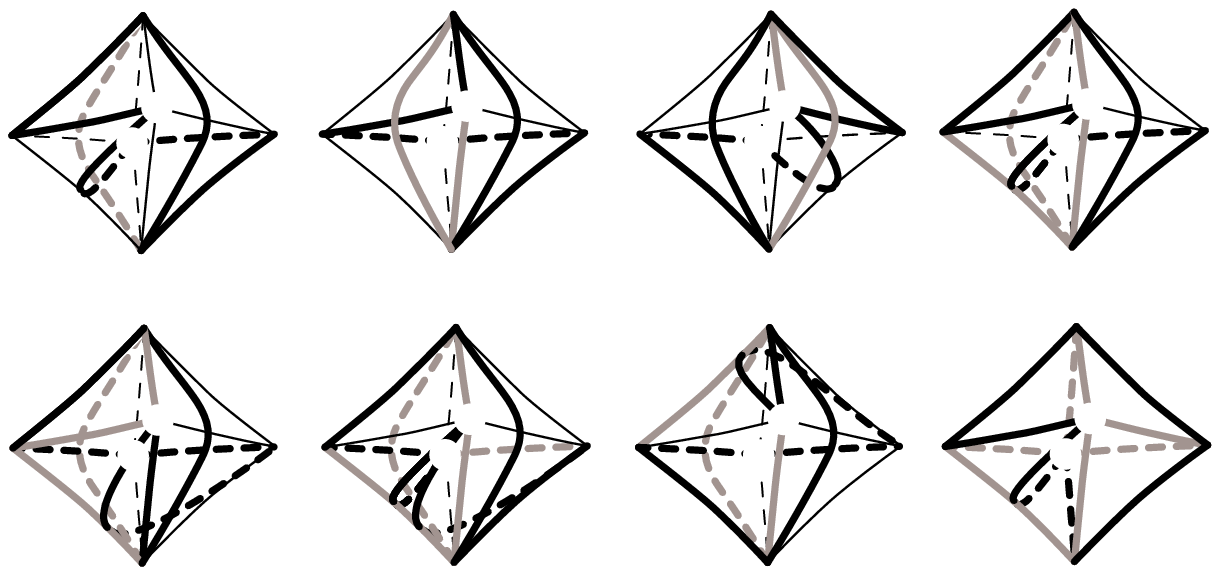}} 
\end{center}
\caption{Applications of (i) or (ii) for
$\la_2=\lb_{2,3}^5,\lb_{2,3}^6$ and $\lb_{2,3}^{6,5}$; of
Theorem~\ref{tri}, (ii) for $\la_2=\lb_{2,4}^5,\lb_{2,4}^6,
\lb_{2,4}^{6,5}$ and $\lb_{1,3}^5$; and of Corollary~\ref{angle corollary} for
$\la_2=\lb_{1,4}^5$}
\label{Fig:arcs_2_1} 
\end{figure}

For $m=1$. If $\{j,k\}\neq\{1,2\}$ or $\{j,k\}\neq\{5,6\}$ then
$l(\ga_1)\leq l(\ga_i),~l(\ga_i)\leq l(\la_1)$ are hypotheses, or
preceding propositions, for some $i\in\{2,3,4\}$. If $\{j,k\}=\{1,2\}$
then, by inspection, $\la_1=\lb_{1,2}^5$ we can again apply argument
(i). By inspection there is no such $\la_1$ for
$\{j,k\}=\{5,6\}$. This completes $P_{m,n,1}$ for $n\in\{5,6\}$.

We now give the arguments for: $\la_m,\ga_m$ share endpoints and
$n>m+1$. The arc set $\lG:=\la_m\cup\ga_m$ divides $\cO$ into two
components. Either: (i) $\lG$ divides one cone point ($c$) from three;
or (ii) $\lG$ divides two cone points from two. For (i) we let
$\cO_c,~\cO'_c$ denote the components of $\cO\sm\lG$ so that
$c\in\cO_c$ and we let $\la'_m$ (respectively $\la''_m$) denote the arc
between $\omega_m,c$ (respectively between $\omega_{m+1},c$) in $\cO_c$.

First $m=4$, (i), $n=6$. None of $\ga_1,\ga_2,\ga_3$ crosses 
$\lG=\la_4\cup\ga_4$, so $\cC_3=\ga_1,\ga_2,\ga_3$ lies in one or other
component of $\cO\sm\lG$. Now $\cC_3$ contains three cone points
disjoint from $\lG$, so $\cC_3\sub\cO'_c$. So $c=\omega_6$ and
$\cC'_4=\ga_1,\ga_2,\ga_3,\la'_4$ is a chain. We observe -- see
Figure~\ref{Fig:arcs_4,3} -- that $\la'_4=\lb_{4,6}$ and hence
$l(\ga_4)\leq l(\la'_4)$ is a hypothesis. By Lemma~\ref{c&p-2}(i):
$2l(\la'_4)<l(\ga_4)+l(\la_4)$ and so $l(\ga_4)\leq
l(\la'_4)<l(\la_4)$.

Second $m=3$, (i), $n\in\{5,6\}$. Neither $\ga_1$ nor $\ga_2$ crosses
$\lG=\la_3\cup\ga_3$, so $\cC_2=\ga_1,\ga_2$ lies in one or other
component of $\cO\sm\lG$. Now $\cC_2$ contains two cone points
disjoint from $\lG$, so $\cC_2\sub\cO'_c,~c=\omega_5$ or $\omega_6$
and $\cC'_3=\ga_1,\ga_2,\la'_3$ is a chain. We observe -- see
Figure~\ref{Fig:arcs_4,3} -- that $\la'_3=\lb_{3,5}$ or
$\la'_3=\lb_{3,6}$ and hence $l(\ga_3)\leq l(\la'_3)$ is
hypothesis. Again, by Lemma~\ref{c&p-2}(i):
$2l(\la'_3)<l(\ga_3)+l(\la_3)$ and so $l(\ga_3)\leq
l(\la'_3)<l(\la_3)$. For (ii) we have that $\la_3=\lb_{3,4}^6$ and
$l(\ga_3)\leq l(\lb_{3,4}^6)$ is a hypothesis.

Next $m=2$, (i), $n\in\{4,5,6\}$. The arc $\ga_1$ does not cross
$\lG=\la_2\cup\ga_2$, so $\ga_1\sub\cO'_c$ and
$c\in\{\omega_4,\omega_5,\omega_6\}$ (respectively $\ga_1\sub\cO_c$ and
$c=\omega_1$). For $n\in\{5,6\}$ -- see Figure~\ref{Fig:arcs_2_1} -- we
have that $\la'_2=\lb_{2,6}$ (respectively $\la''_2=\lb_{1,3}$). For $n=4$ -- %
see Figure~\ref{Fig:arcs_2_2} -- we have that $\la'_2=\lb_{2,4}$ or
$\lb_{2,5}$ (respectively there is no such $\la_2$). So $l(\ga_2)\leq
l(\la'_2)$ (respectively $l(\ga_2)\leq l(\la''_2)$) is a hypothesis. By
Lemma~\ref{c&p-2}(i): $2l(\la'_2)$ or $2l(\la''_2)<l(\ga_2)+l(\la_2)$
and so $l(\ga_2)\leq l(\la'_2)<l(\la_2)$ (respectively $l(\ga_2)\leq
l(\la''_2)<l(\la_2)$).

For (ii), again, $\ga_1$ lies in one component of $\cO\sm\lG$. Let
$\la'''_2$ denote the unique arc disjoint from $\lG$ in this component
of $\cO\sm\lG$. For $n\in\{5,6\}$ -- again see Figure~\ref{Fig:arcs_2_1}%
 -- we have that $\la'''_2=\ga_6$. For $n=4$ -- again see
Figure~\ref{Fig:arcs_2_2} -- we have $\la'''_2=\lb_{1,4}$ or
$\lb_{1,5}$. So $l(\ga_2)\leq l(\la'''_2)$ is a hypothesis. By
Lemma~\ref{c&p-2}(ii): $2l(\la'''_2)<l(\ga_2)+l(\la_2)$ and so
$l(\ga_2)\leq l(\la'''_2)<l(\la_2)$.

Finally, $m=1$, (i), $n\in\{3,\ldots,6\}$. For $n\in\{5,6\}:~\la'_1=\ol{\lb_{2,6}}$ and $l(\ga_2)\leq l(\la'_1)$ is a
hypothesis. For $n\in\{3,4\}:~l(\ga_2)\leq l(\la'_1)$
is a proceeding proposition. Since $l(\ga_1)\leq l(\ga_2)$ is a
hypothesis, we have that $l(\ga_1)\leq l(\ga_2)\leq l(\la'_1)$. By
Lemma~\ref{c&p-2}(i): $2l(\la'_1)<l(\ga_1)+l(\la_1)$ and so
$l(\ga_1)\leq l(\la'_1)<l(\la_1)$. 

For (ii), $n\in\{5,6\}$, there is no such $\la_1$. For $n\in\{3,4\}$,
we let $\la'_3$ denote the unique arc disjoint from $\lG$ in the same
component of $\cO\sm\lG$ as $\ga_2$. Here $\cC'_3=\ga_1,\ga_2,\la'_3$
is a chain and so $l(\ga_3)\leq l(\la'_3)$ is a proceeding
proposition. Since $l(\ga_1)\leq l(\ga_3)$ is a hypothesis, we have
that $l(\ga_1)\leq l(\ga_3)\leq l(\la'_3)$ . By Lemma~\ref{c&p-2}(ii):
$2l(\la'_3)<l(\ga_1)+l(\la_1)$ and so $l(\ga_1)\leq
l(\la'_3)<l(\la_1)$.

Now suppose $n\in\{1,\ldots,4\},~\la_m$ crosses $\cC_4$. 

\begin{lemma}\label{no.>1}  

Suppose that either $X(\la_m,\ga_n)>1$ or $\la_m,\ga_n$ share an
endpoint. Then there exist arcs $\la'_m,\ga'_n$ between the same
respective endpoints as $\la_m,\ga_n$ such that $l(\la'_m)<l(\la_m)$
or $l(\ga'_n)<l(\ga_n);~
X(\la'_m,\ga_n),X(\ga'_n,\ga_n)<X(\la_m,\ga_n)$; and
$X(\la'_m,\ga_i)=X(\ga'_n,\ga_i)=0$ for $i\leq n-1$. In particular
$\cC'_m=\ga_1,\ldots,\ga_{m-1},\la'_m,
\cC''_n=\ga_1,\ldots,\ga_{n-1},\ga'_n$ are both chains.

\end{lemma}

\begin{proof} This result is essentially Proposition~3.1 in
\cite{Gri2}, with additional observations upon the number of crossing
points. However, upon going through the proof, these observations
become clear.\end{proof}

The following argument gives $P_{m,n,p>1}$: it uses induction on
$p$, the first induction step being the set of propositions that precede
$P_{m,n,p>1}$.

Let $X(\la_m,\ga_n)=p>1$ and so by Lemma~\ref{no.>1} there
exist arcs $\la'_m,\ga'_n$ as stated. Let
$p'=X(\la'_m,\ga_n)<p, ~p''=X(\ga'_m,\ga_n)<p$. We note
that $l(\ga_m)\leq l(\la'_m)$ is either: $P_{m,n,p'>1}$ if
$p'>1$; or a preceding proposition if $p'\leq 1$. Likewise,
$l(\ga_n)\leq l(\ga'_n)$ is either: $P_{m,n,p''>1}$ if
$n=m$ and $p''>1$; or a preceding proposition if $n>m$ or
$p''\leq 1$. Since $l(\la'_m)<l(\la_m)$ or
$l(\ga'_n)<l(\ga_n)$ it follows, by induction on $p$, that
$l(\ga_m)\leq l(\la'_m)<l(\la_m)$.

So, for the rest of the proof, we may suppose that $X(\la_m,\ga_n)=1$.

\begin{lemma}\label{no.=1}  

Suppose that $\la_m,\ga_n$ have distinct endpoints and that
$k>n+1$. Then there exist arcs $\la'_m,\ga'_n$ between
$\omega_j,\omega_{n+1}$ and $\omega_n,\omega_k$ such that
$l(\la'_m)<l(\la_m)$ or $l(\ga'_n)<l(\ga_n)$ and
$X(\la'_m,\ga_i)=X(\ga'_n,\ga_i)=0$ for $i\leq n$. In particular
$\cC'_m=\ga_1,\ldots,\ga_{m-1},\la'_m,
\cC''_n=\ga_1,\ldots,\ga_{n-1},\ga'_n$ are both chains.

\end{lemma}

\begin{proof} This is essentially Lemma~3.3 in \cite{Gri2}, again with
additional observations upon the number of crossing points. Again, these
observations are clear.\end{proof}

We now give two general arguments using these two lemmas.

Suppose: (1)\qua $\la_m,\ga_n$ share an endpoint. Again we can apply
Lemma~\ref{no.>1}: there exist arcs $\la'_m,\ga'_n$ as stated. In
particular $X(\la'_m,\ga_i)=X(\ga'_n,\ga_i)=0$ for $i\leq n$. So
$l(\ga_m)\leq l(\la'_m),~l(\ga_n)\leq l(\ga'_n)$ are both preceding
propositions. Since $l(\la'_m)<l(\la_m)$ or $l(\ga'_n)<l(\ga_n)$, it
follows that $l(\ga_m)\leq l(\la'_m)<l(\la_m)$.

Suppose: (2)\qua $\la_m,\ga_n$ have distinct endpoints and $k>n+1$. By
Lemma~\ref{no.=1} there exist arcs $\la'_m,\ga'_n$ as stated. Again
$l(\ga_m)\leq l(\la'_m),l(\ga_n)\leq l(\ga'_n)$ are both preceding
propositions. As $l(\la'_m)<l(\la_m)$ or $l(\ga'_n)<l(\ga_n)$, we have
that $l(\ga_m)\leq l(\la'_m)<l(\la_m)$.

For $m=4:j=4,k\in\{5,6\}$ and $n=4:\la_4,\ga_4$ share the
endpoint $\omega_4$ (1).

For $m=3:j=3,k\in\{4,5,6\}$. For $n=4$ if $k\in\{4,5\}$ then
$\la_3,\ga_4$ share the endpoint $\omega_k$ (1); if $k=6$ then
$\la_3,\ga_4$ have distinct endpoints and $k>n+1$ (2). For
$n=3:\la_3,\ga_3$ share the endpoint $\omega_3$ (1).

For $m=2:j\in\{1,2\},k\in\{3,\ldots,6\}$. For $n=4$ if $k=3$ then, by
inspection, $\la_2$ is one of
$\lb_{2,3}^4,\lb_{2,3}^{4,5,6},\lb_{2,3}^{5,4},\lb_{2,3}^{6,4},
\lb_{2,3}^{6,5,4}$, and we can apply argument (i) or (ii), or is one of
$\lb_{1,3}^4,\lb_{1,3}^{4,5,6},\lb_{1,3}^{5,4}$, and we apply
Theorem~\ref{tri} in conjunction with argument (ii) -- see Figure~5. If
$k\in\{4,5\}$ (1); if $k=6$ (2). For $n=3$ if $k\in\{3,4\}$ (1); if
$k\in\{5,6\}$ (2). For $n=2$ if $k=3$ (1); if $k\in\{4,5,6\}$ (2).

\begin{figure}[htbp] 
\begin{center} 
\leavevmode 
\SetLabels 
(-0.0125*0.775)\mbox{\scriptsize $\omega_1$}\\
(0.115*1.01)\mbox{\scriptsize $\omega_2$}\\
(0.115*0.54)\mbox{\scriptsize $\omega_3$}\\
(0.1275*0.825)\mbox{\scriptsize $\omega_4$}\\
(0.23*0.745)\mbox{\scriptsize $\omega_5$}\\
(0.1075*0.7575)\mbox{\scriptsize $\omega_6$}\\
(0.385*0.82)\mbox{\scriptsize $\omega_1$}\\
(0.3725*1.01)\mbox{\scriptsize $\omega_2$}\\
(0.3725*0.54)\mbox{\scriptsize $\omega_3$}\\
(0.4875*0.74)\mbox{\scriptsize $\omega_4$}\\
(0.26*0.8)\mbox{\scriptsize $\omega_5$}\\
(0.365*0.7575)\mbox{\scriptsize $\omega_6$}\\
(0.51*0.8)\mbox{\scriptsize $\omega_1$}\\
(0.6325*1.01)\mbox{\scriptsize $\omega_2$}\\
(0.6325*0.54)\mbox{\scriptsize $\omega_3$}\\
(0.745*0.74)\mbox{\scriptsize $\omega_4$}\\
(0.6425*0.82)\mbox{\scriptsize $\omega_5$}\\
(0.6225*0.7575)\mbox{\scriptsize $\omega_6$}\\
(0.8975*0.82)\mbox{\scriptsize $\omega_1$}\\
(0.8875*1.01)\mbox{\scriptsize $\omega_2$}\\
(0.8875*0.54)\mbox{\scriptsize $\omega_3$}\\
(0.7725*0.8)\mbox{\scriptsize $\omega_4$}\\
(0.8775*0.7575)\mbox{\scriptsize $\omega_5$}\\
(1.01*0.775)\mbox{\scriptsize $\omega_6$}\\
(0.125*0.25)\mbox{\scriptsize $\omega_1$}\\
(0.115*0.44)\mbox{\scriptsize $\omega_2$}\\
(0.115*-0.03)\mbox{\scriptsize $\omega_3$}\\
(0.105*0.185)\mbox{\scriptsize $\omega_4$}\\
(-0.0125*0.20)\mbox{\scriptsize $\omega_5$}\\
(0.23*0.18)\mbox{\scriptsize $\omega_6$}\\
(0.3725*-0.03)\mbox{\scriptsize $\omega_1$}\\
(0.3725*0.44)\mbox{\scriptsize $\omega_2$}\\
(0.26*0.24)\mbox{\scriptsize $\omega_3$}\\
(0.38375*0.25)\mbox{\scriptsize $\omega_4$}\\
(0.36375*0.185)\mbox{\scriptsize $\omega_5$}\\
(0.4875*0.175)\mbox{\scriptsize $\omega_6$}\\
(0.6325*0.44)\mbox{\scriptsize $\omega_1$}\\
(0.6325*-0.03)\mbox{\scriptsize $\omega_2$}\\
(0.5125*0.24)\mbox{\scriptsize $\omega_3$}\\
(0.745*0.175)\mbox{\scriptsize $\omega_4$}\\
(0.64375*0.25)\mbox{\scriptsize $\omega_5$}\\
(0.6225*0.185)\mbox{\scriptsize $\omega_6$}\\
(0.885*-0.03)\mbox{\scriptsize $\omega_1$}\\
(0.885*0.44)\mbox{\scriptsize $\omega_2$}\\
(0.7725*0.24)\mbox{\scriptsize $\omega_3$}\\
(0.875*0.185)\mbox{\scriptsize $\omega_4$}\\
(0.895*0.25)\mbox{\scriptsize $\omega_5$}\\
(1.01*0.2)\mbox{\scriptsize $\omega_6$}\\
\endSetLabels 
\AffixLabels{\epsffile{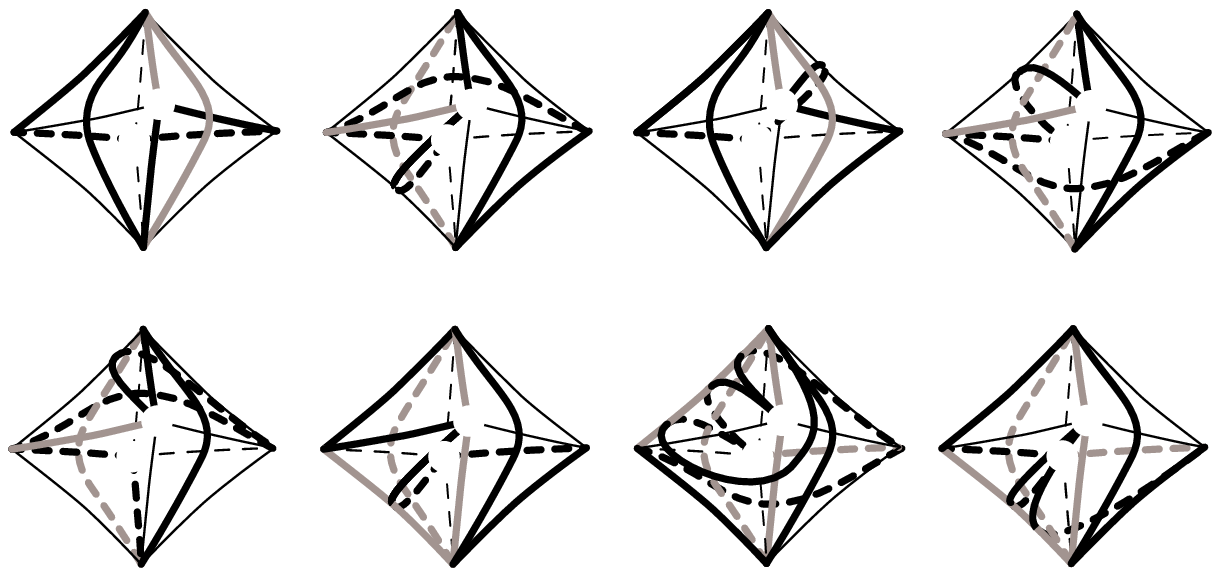}} 
\end{center}
\caption{For
$\la_2=\lb_{2,3}^4,\lb_{2,3}^{4,5,6},\lb_{2,3}^{5,4}, \lb_{2,3}^{6,4}$
and $\lb_{2,3}^{6,5,4}$ applications of (i) or (ii); and for
$\la_2=\lb_{1,3}^4,\lb_{1,3}^{4,5,6}$ and $\lb_{1,3}^{5,4}$
applications of Theorem~\ref{tri}, (ii)}
\label{Fig:arcs_2_2} 
\end{figure}

Finally $m=1$. Suppose $n=4$. If $\{j,k\}\neq\{1,2\}$ or
$\{j,k\}\neq\{5,6\}$ then $l(\ga_1)\leq l(\ga_i),l(\ga_i)\leq
l(\la_1)$ are both preceding propositions for some $i\in\{2,3,4\}$. If
$\{j,k\}=\{1,2\}$ we can apply (i) or (ii). There
is no such $\la_1$ for $\{j,k\}=\{5,6\}$.
 
Now suppose $n=3$. If $\{j,k\}\neq\{1,2\}$ or
$\{j,k\}\not\sub\{4,5,6\}$ then $l(\ga_1)\leq l(\ga_i),l(\ga_i)\leq
l(\la_1)$ are both preceding propositions for some
$i\in\{2,3\}$. Again, if $\{j,k\}=\{1,2\}$ we can apply (i) or
(ii). For $\{j,k\}\sub\{4,5,6\}$ either $j=4$ (1) or $j=5$ (2).

Now suppose $n=2$. If $\{j,k\}\neq\{1,2\}$ or
$\{j,k\}\not\sub\{3,\ldots,6\}$ (ie $j\in\{1,2\},k\in\{3,\ldots,6\}$)
then $l(\ga_1)\leq l(\ga_2),l(\ga_2)\leq l(\la_1)$ are both preceding
propositions. For $\{j,k\}=\{1,2\}$ (1). For $\{j,k\}\sub\{3,\ldots,6\}$
either $j=3$ (1); or $j\in\{4,5,6\}$ (2).

Finally $n=1$. Either $j$ or $k\in\{1,2\}$ (1); or
$\{j,k\}\sub\{3,\ldots,6\}$ (2).\end{proof}

\begin{proof}[Proof of Theorem~\ref{angle characterisation}] As
$l(\lk_{3,0})\leq l(\lk_{0,5}),l(\lk_{2,3})\leq l(\lk_{2,5}),l(\lk_{0,1})\leq
l(\lk_{0,4})$, by Corollary~\ref{obvious}, we have that $l(\lk_{1,2})\geq
l(\lk_{2,4})$. Likewise, since $l(\lk_{3,0})\leq l(\lk_{0,4}),l(\lk_{2,3})\leq
l(\lk_{2,4}),l(\lk_{0,1})\leq l(\lk_{0,5})$ we have that $l(\lk_{1,2})\geq
l(\lk_{2,5})$. That is $l(\lk_{1,2})\geq l(\lk_{2,l})$.

The arc set $K$ divides $\cO$ into eight triangles. We label these as
follows: let $t_k$ (respectively $T_k$) denote the triangle with one edge
$\lk_{k,k+1}$ and one vertex $c_4$ (respectively $c_5$). We shall use $\angle
c_l t_k$ to denote the angle at the $c_l$--vertex of $t_k$, et
cetera. Cut $\cO$ open along
$\lk_{3,0}\cup\lk_{0,1}\cup\lk_{1,4}\cup\lk_{1,2}\cup\lk_{1,5}$ to
obtain a domain $\lO$.

We show that $l(\lk_{2,3})\leq l(\lk_{2,l}),l(\lk_{3,0})\leq
l(\lk_{1,2})\leq \{l(\lk_{0,l}),l(\lk_{1,l})\},l(\lk_{0,1})\leq
l(\lk_{0,l})$ implies that ${\rm min}_l~l(\lk_{3,l})\leq l(\lk_{0,1})$ with
equality {\em if and only if} $\cO$ is the octahedral orbifold. First
we show that: $\angle c_2 t_2\leq\angle c_4 t_0~or~\angle c_2
T_2\leq\angle c_5 T_0$.

Now $l(\lk_{1,2})\leq l(\lk_{1,l})$, $l(\lk_{3,0})\leq
l(\lk_{0,l}$ so 
$\angle c_2 t_1\geq\angle c_4 t_1$, $\angle c_2 T_1\geq\angle c_5 T_1$,
$\angle c_3 t_3\geq\angle c_4 t_3$, $\angle c_3 T_3\geq\angle c_5 T_3$,  
which imply 
$$\begin{aligned}\angle c_2 t_1+\angle c_2 T_1 + \angle c_3 t_3 + \angle c_3 T_3&\geq\angle c_4 t_1 + \angle
c_5 T_1 + \angle c_4 t_3 + \angle c_5 T_3\\  
\Lr(\pi-\angle c_2 t_1-\angle
c_2 T_1) + (\pi-\angle c_3 t_3 -& \angle c_3 T_3)\\
\leq(\pi-\angle c_4 t_1 -&
\angle c_4 t_3) + (\pi-\angle c_5 T_1 - \angle c_5 T_3)\\  
\Lr(\angle c_2
t_2+\angle c_2 T_2) + (\angle c_3 t_2 + \angle c_3 T_2)&\leq(\angle c_4 t_2 +
\angle c_4 t_0) + (\angle c_5 T_2 + \angle c_5 T_0)\end{aligned}$$ 
and $l(\lk_{2,3})\leq l(\lk_{2,l})$ so $\angle c_3 t_2\geq \angle c_4
t_2,\angle c_3 T_2\geq \angle c_5 T_2$ 
$\Ra\angle c_2 t_2+\angle c_2 T_2\leq \angle c_4 t_0 + \angle c_5 T_0$
$\Ra\angle c_2 t_2\leq\angle c_4 t_0$ or $\angle c_2 T_2\leq\angle c_5 T_0.$

\begin{figure}[htbp] 
\begin{center} 
\leavevmode 
\SetLabels 
(0.15*0.5)\mbox{\scriptsize $c_4$}\\
(0.18*-0.01)\mbox{\scriptsize $\ol{c_1}$}\\
(0.51*0.1)\mbox{\scriptsize $c_2$}\\
(0.51*0.89)\mbox{\scriptsize $c_3$}\\
(0.2*1.0)\mbox{\scriptsize $c_0$}\\
(0.0*0.8)\mbox{\scriptsize $c_1$}\\
(0.875*0.5)\mbox{\scriptsize $c_5$}\\
(0.85*-0.01)\mbox{\scriptsize $\ol{c'_1}$}\\
(0.83*1.0)\mbox{\scriptsize $c_0$}\\
(1.01*0.8)\mbox{\scriptsize $c'_1$}\\
(0.15*0.76)\mbox{\scriptsize $t_0$}\\
(0.3*0.21)\mbox{\scriptsize $t_1$}\\
(0.4*0.5)\mbox{\scriptsize $t_2$}\\
(0.3*0.77)\mbox{\scriptsize $t_3$}\\
(0.87*0.76)\mbox{\scriptsize $T_0$}\\
(0.7*0.21)\mbox{\scriptsize $T_1$}\\
(0.62*0.5)\mbox{\scriptsize $T_2$}\\
(0.7*0.77)\mbox{\scriptsize $T_3$}\\
\endSetLabels 
\AffixLabels{ \epsffile{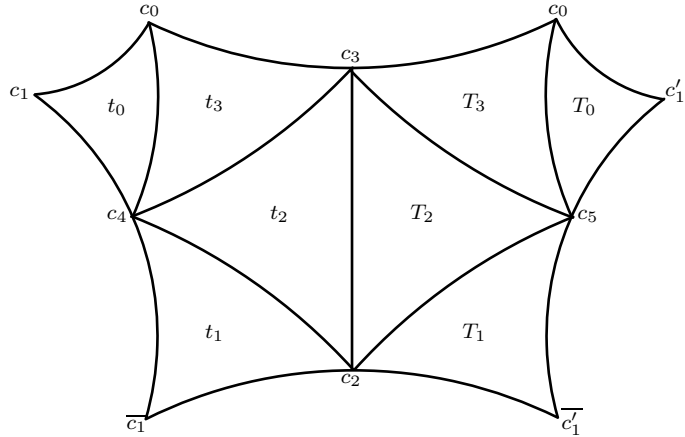}} 
\end{center}
\caption{The triangles $t_k,T_k$ in the domain $\lO$}
\label{Fig:t1_T1} 
\end{figure}
 
Up to relabelling, we may suppose that $\angle c_2 t_2\leq\angle c_4 t_0$. We
now show that $l(\lk_{3,4})\leq l(\lk_{0,1})$. There are two arguments. Firstly
we show that if $\angle c_3 t_2\geq \pi -\theta$ then
$l(\lk_{0,4})<l(\lk_{3,0})$ -- contradicting a hypothesis. So $\angle c_3
t_2<\pi -\theta$ and we then show that $l(\lk_{3,4})\leq l(\lk_{0,1})$. The
angle $\theta$ is given as follows. Let $\cI_2$ be an isoceles triangle with
vertices $v_2,v_3,v_4$ and edges $\lv_{2,3},\lv_{2,4},\lv_{3,4}$ such that
$l(\lv_{2,3})=l(\lv_{2,4})=l(\lk_{2,4})$ and $\angle v_2\cI_2=\angle c_2 t_2$.
Then $\theta=\angle v_3\cI_2=\angle v_4\cI_2$. 

Let $C_2,C_4$ denote circles of radius $l(\lk_{2,4})$ about $c_2,c_4$
respectively. As in Figure~\ref{Fig:triangles} $c_3$ must lie inside
$C_2$ since $l(\lk_{2,3})\leq l(\lk_{2,4})$. Likewise $c_0$ must lie
outside $C_4$ since $l(\lk_{0,4})\geq l(\lk_{1,2})\geq
l(\lk_{2,4})$. Similarly $c_1$ must lie outside $C_4$ since
$l(\lk_{1,4})\geq l(\lk_{1,2})\geq l(\lk_{2,4})$.  Moreover since the
angle sum at any cone point is $\pi:~\angle c_3 t_2+\angle c_3
t_3<\pi$. In Figure~6 we have also constructed the point $x$ as the
intersection of the radius through $\lk_{2,3}$ and $\cC_4$. Let $t_x$
denote the triangle spanning $x,c_3,c_4$.

Now $\angle c_3 t_2\geq\pi - \theta$ is equivalent to $\angle c_3
t_x\leq\theta$. It follows that $\angle c_4 t_x\geq\angle c_3 t_x$. By
inspection $\angle c_4 t_3>\angle c_4 t_x$ and $\angle c_3 t_x>\angle c_3 t_3$.
So $\angle c_4 t_3>\angle c_4 t_x\geq\angle c_3 t_x>\angle c_3 t_3$ or
equivalently $l(\lk_{0,4})<l(\lk_{0,3})$.

So $\angle c_3 t_2<\pi-\theta$ and we will compare $t_2,t_0$. Firstly, $\angle
c_3 t_2<\pi-\theta$ implies that $l(\lk_{3,4})\leq l(\lv_{3,4})$. (Recall that
$\lv_{3,4}$ is an edge of $\cI_2$.) Let $\cI_0$ be an isoceles triangle with
vertices $v_0,v_1,v_4$ and edges $\lv_{0,1},\lv_{1,4},\lv_{0,4}$ such that
$l(\lv_{1,4})=l(\lv_{0,4})=l(\lk_{2,4})$ and $\angle v_4\cI_0=\angle c_4 t_0$.
Since $l(\lk_{0,4}),l(\lk_{1,4})\geq l(\lk_{1,2})\geq l(\lk_{2,4})$ we then
observe that $l(\lk_{0,1})\geq l(\lv_{0,1})$. As $\angle c_2 t_2\leq\angle c_4
t_0$ we have that $l(\lv_{3,4})\leq l(\lv_{0,1})$. Therefore $l(\lk_{0,1})\geq
l(\lv_{0,1})\geq l(\lv_{3,4})\geq l(\lk_{3,4})$.

We have equality {\em if and only if} $\angle c_2 t_2=\angle c_4 t_0$ and
$l(\lk_{2,3})=l(\lk_{2,4})=l(\lk_{0,4})=l(\lk_{1,4})$. From above $\angle c_2
t_2=\angle c_4 t_0$  {\em if and only if} $l(\lk_{1,2})=l(\lk_{1,l}),
l(\lk_{3,0})=l(\lk_{0,l})$ and $l(\lk_{2,3})=l(\lk_{2,l})$. So we have that 
$l(\lk_{0,1})=l(\lk_{3,4})$ and
$l(\lk_{1,2})=l(\lk_{2,3})=l(\lk_{3,0})=l(\lk_{0,l})=l(\lk_{1,l})=l(\lk_{2,l})$.

That is: $t_1,T_1$ are isometric equilateral triangles and 
$t_0,T_0,t_2,t_3$ (respectively $T_2,T_3$) are isometric isoceles triangles. By
considering angle sums at $c_4,c_5:\angle c_4 t_2=\angle c_4 t_3=\angle c_5
T_2=\angle c_5 T_3$. So: $t_1,T_1$ are isometric equilateral triangles
and $t_0,T_0,t_2,t_3,T_2,T_3$ are isometric isoceles triangles. By the
angle sum at $c_3:\angle c_3 t_2=\angle c_3 t_3=\angle c_3 T_2=\angle c_3
T_3=\pi/4$ and so $\angle c_0 t_0=\angle c_1 t_0=\angle c_0 T_0=\angle c_1
T_0=\pi/4$. Again, by considering angle sums at $c_0,c_1$ all the angles are
$\pi/4$, all of the edges are of equal length. So $\cO$ is the octahedral
orbifold.\end{proof}

\begin{figure}[htbp] 
\begin{center} 
\leavevmode 
\SetLabels 
(0.045*0.81)\mbox{\scriptsize $c_0$}\\
(0.38*0.85)\mbox{\scriptsize $c_2$}\\
(0.27*0.8)\mbox{\scriptsize $c_3$}\\
(0.2175*0.475)\mbox{\scriptsize $c_4$}\\
(0.08*0.84)\mbox{\scriptsize $x$}\\
(0.62*0.81)\mbox{\scriptsize $c_0$}\\
(0.59*0.55)\mbox{\scriptsize $c_1$}\\
(0.9525*0.85)\mbox{\scriptsize $c_2$}\\
(0.78*0.8)\mbox{\scriptsize $c_3$}\\
(0.8*0.46)\mbox{\scriptsize $c_4$}\\
(0.67*0.62)\mbox{\scriptsize $t_0$}\\
(0.83*0.68)\mbox{\scriptsize $t_2$}\\
(0.73*0.7)\mbox{\scriptsize $t_3$}\\
\endSetLabels 
\AffixLabels{ \epsffile{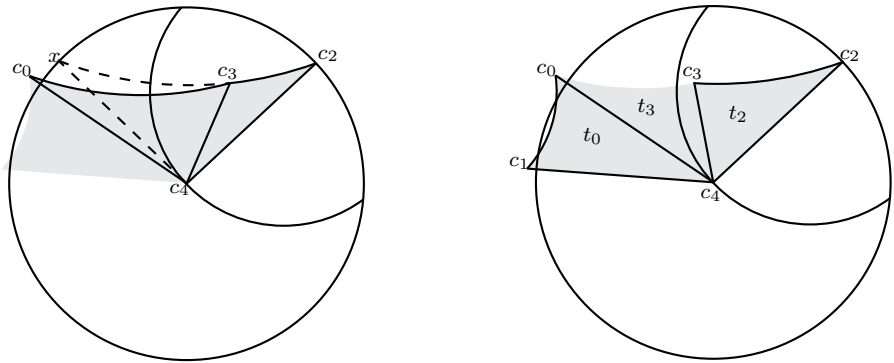}} 
\end{center}
\caption{Arguments for $\angle c_3 t_2\geq\pi-\theta$ and
for $\angle c_3 t_2<\pi-\theta$}
\label{Fig:triangles} 
\end{figure}
\np

\Addresses\recd

\end{document}